\documentclass[10pt,a4paper,draft]{article}
\usepackage{amsmath}
\usepackage{amsfonts}
\usepackage{amssymb}
\usepackage{amscd}
\makeatletter\@addtoreset{equation}{section}\makeatother

\newtheorem{theorem}{Theorem}[section]

\newtheorem{lemma}{Lemma}
\newenvironment{corollary1}{\par\noindent\ignorespaces{\bf
Corollary 1}\,\sl}

\def\C{{\mathbb C}}

\newtheorem{definition}{Definition}

\begin{document}
\title{{\bf
Non-unitary set-theoretical solutions to the Quantum Yang-Baxter
Equation}}
\author{ {\Large Alexander Soloviev}\\
MIT \\
Department of Mathematics\\
Cambridge, MA 02139,
      USA
\date{}
}
\maketitle
\begin{abstract}
We develop a theory of non-unitary set-theoretical solutions to the
Quantum Yang-Baxter equation. Our results generalize those obtained by
Etingof, Schedler and the author in \cite{ESS}. We remark that some of
our constructions are similar to constructions obtained by Lu, Yan and Zhu
in \cite{LYZ}.
\end{abstract}
\section{Introduction}
In this paper we study set-theoretical solutions to the Quantum
Yang-Baxter equation, i.e. permutations $R:X\times X\to X\times X$
with $X$ being a non-empty set such that
$$
R^{12}R^{13}R^{23}=R^{23}R^{13}R^{12} \ in\ Aut(X\times X\times X).
$$
In the above $R^{12}$, $R^{13}$, $R^{23}\in  Aut(X\times X\times X)$
stand for $R$ acting in 1,2; 1,3; and 2,3 components of 
$X\times X\times X$ correspondingly. The idea to consider
set-theoretical solutions first appeared in \cite{Dr}. Later on,
Etingof, Schedler and the author \cite{ESS} studied set-theoretical
solutions to the
Quantum Yang-Baxter equation which satisfied additional properties of
unitarity and nondegeneracy (crossing symmetry). In particular,
\cite{ESS} contained the
classification of nondegenerate unitary set-theoretical solutions to QYBE
in group theoretical terms as well as numerous classes of examples of such
solutions. Subsequently, Lu, Yan and Zhu in \cite{LYZ} showed that many
of the constructions from \cite{ESS} hold in a more general case of
nondegenerate but not necessarily unitary set-theoretical solutions to
QYBE. 

Following \cite{ESS} and \cite{LYZ}, we develop a theory of nondegenerate
set-theoretical solutions to QYBE. Particularly, we show that the 
unitarity condition that was used in \cite{ESS} for group theoretical
characterization
of unitary nondegenerate set-theoretical solutions to QYBE can be dropped.
We give a group theoretical characterization of general set-theoretical
nondegenerate solutions to QYBE in Theorem \ref{1-1}.  
We also introduce injective solutions, study their properties, and show
that there is a combinatorial criterion  (Theorem  \ref{criterion}) 
describing the class of injective nondegenerate solutions  to
QYBE. Injectivity property, which is a generalization of involutivity, is
important for studying affine solutions. In
particular the classification of unitary affine solutions given in
\cite{ESS} can be generalized to include injective solutions.  

It was shown in \cite{ESS}
that the structure group of a  nondegenerate unitary set-theoretical
solution to QYBE on a set with N elements always has an abelian subgroup
of finite index and of rank N. We compute (Theorem \ref{rank}) the rank of
a
finite index abelian subgroup of
the structure group for an arbitrary nondegenerate finite solution and
show that
this rank never exceeds N, with the equality taking place only in the
unitary case. 

In the second part of the paper we discuss the applications
of the developed theory to examples. In particular, we classify affine
solutions on an abelian group. It is proved that injective affine
solutions are obtained from the representations of the algebra generated
by invertible elements $p$, $q$, $z$ subject to $pq=qp$ and
$z^2-z(p+q)+pq=0$. 

{\bf Acknowledgments.}
The author is thankful to his advisor Pavel Etingof for a valuable
exchange of ideas and guidance throughout the work on the paper, and to
Lu, Yan, Zhu for making their work available before publication.  

\section{Structure groups}
\subsection{Construction of the structure groups}

Let $X$ be a nonempty set and $S:X\times X\to X\times X$ a
bijective map. We call a pair $(X,\ S)$ a braided set if the following
braiding condition holds in $X\times X\times X$:

\begin{equation}
\label{braid}
S_1S_2S_1=S_2S_1S_2,
\end{equation}
where $S_1=S\times id$, $S_2=id\times S$.

{\noindent \underline{Remark.}} Consider the map 
$R:X\times X\to X\times X$ given by $R=\sigma S$, where
$\sigma(x,y)=(y,x)$ for $x,\ y\in X$. Then $(X, S)$ is a braided set
if and only if $R$ satisfies the Quantum Yang-Baxter equation. 

It is useful to associate with a braided set $(X,S)$ two
groups $G_X$ and $A_X$. 
\begin{definition}
Define the group $G_X$ as the group generated
by the elements of $X$ subject to the relations $xy=y_1x_1$ 
if $S(x,y)=(y_1,x_1)$, where  $x, y \in
X$. We call $G_X$ the structure group of the braided set $(X,\ S)$.
\end{definition}
\begin{definition}
Define the group $A_X$
as the group  generated by elements of $X$ subject to relations 
$x_1\bullet y=y_2\bullet x_1$,  where $x, y \in X$ and $x_1,\ y_2\in
X$ are defined out of relations $S(x,y)=(y_1,x_1),\
S(y_1,x_1)=(x_2,y_2)$. We
call $A_X$ the derived structure group of the braided set $(X,\ S)$. 
\end{definition}
We introduce the
maps $g:X\times X\to X$ and $f:X\times X\to X$ as
components of $S$, i.e. for $x,\ y\in X$
 $$S(x,y)=(g_x(y),f_y(x)).$$
\begin{definition}
(i) We call a set $(X,S)$ nondegenerate if $g_x(y)$
is a bijective function of $y$ for fixed $x$ and $f_y(x)$ is a bijective
function of
$x$ for fixed $y$.
(ii) We call a  set $(X,S)$ involutive if $S^2=id_{X^2}$. A
braided set which is involutive will be called symmetric.    
\end{definition}
In particular, for a symmetric set $(X,S)$ we see that the group $A_X$
is the free abelian group generated by elements of $X$. Note that the
properties of involutivity and nondegeneracy are equivalent to
corresponding properties of unitarity and crossing symmetry for the map
$R=\sigma S$ \cite{ESS}. 

Recall that the braid group $B_n$ for $n\geq 2$ is generated by elements
$b_i,\ 1\leq
i\leq n-1$, with defining relations 
$$
b_ib_j=b_jb_i,\ |i-j|>1,\ b_ib_{i+1}b_i=b_{i+1}b_ib_{i+1},
$$
and recall that the symmetric group $S_n$ is the quotient of $B_n$ by the
relations $b_i^2=1$. The following obvious proposition explains
our terminology. Let $S^{ii+1}_n:X^n\to X^n$ be defined as 
$S^{ii+1}_n=id_{X^{i-1}}\times S\times id_{X^{n-i-1}}$.
\begin{theorem}[\cite{ESS}]
\label{twisted}
(i) The assignment $b_i\to S^{ii+1}_n$ extends to an action of $B_n$ on
$X^n$ ($n\geq 3)$ if and only if $(X,\ S)$ is a braided set. 

(ii) The assignment $b_i\to S^{ii+1}_n$ extends to an action of $S_n$ on
$X^n$ ($n\geq 3$) if and only if $(X,\ S)$ is a symmetric set. 
\end{theorem}
\begin{definition}
The action of Theorem \ref{twisted} is called the twisted action of 
$B_n$ (or $S_n$) given by $S$.
\end{definition}
The following
result shows how a nondegenerate braided set $(X,\ S)$ gives rise to two
actions of the structure group $G_X$ on the set $X$.
\begin{theorem}(\cite{ESS})
Suppose that $(X,S)$ is nondegenerate.
Then $(X,S)$ is a braided set
if and only if the following conditions are simultaneously satisfied:
  
(i)  the assignment $x\to f_x$ is a right action of
$G_X$ on $X$;

(ii) the assignment $x\to g_x$ is a left action
of $G_X$ on $X$;

(iii) the linking relation
$$
f_{g_{f_y(x)}(z)}(g_x(y))=g_{f_{g_y(z)}(x)}(f_z(y))
$$
holds.
\label{linking}
\end{theorem}
{\bf Proof:} 
Conditions (i)-(iii) are exactly the three components of the braid
relation (\ref{braid}). $\blacksquare$

In light of the above proposition it makes sense to introduce the
notations $x\circ y=g_x(y)$ and $y*x=f_y^{-1}(x)$ for $x, y\in X$. Then
if $(X, S)$ is a nondegenerate braided set, we can extend $*$ and $\circ$
to left actions of $G_X$ on $X$. We will denote the actions of an
element $g\in G_X$ on an element $x\in X$ by $g\circ x$ and $g*x$
correspondingly.

From  now on we always assume $(X,\ S)$ to be
a nondegenerate braided set. Sometimes we refer to nondegenerate braided
sets as to "solutions", meaning that $S$ is a solution to braid equation
(\ref{braid}).

Define  $\phi:X\times X\to X$ by 
\begin{equation}
\label{phi}
\phi(y,x)=x^{-1}*((y*x)\circ y)
\end{equation}
and
$S':X\times X\to X\times X$ as $S'(x,y)=(\phi(y,x),x)$.
\begin{theorem}
\label{derived}
(i) $\phi(y,z)$ is $G_X$-invariant w.r.t. *-action, i.e. for $g\in
G_X$, $g*\phi(y,z)=\phi(g*y,g*z)$. 

(ii) $(X, S')$ is a nondegenerate braided set. We call this set the
derived braided set or the derived solution.

(iii) The structure group of the derived solution is the derived structure
group.

(iv) For each integer $n\geq 2$ there exists a bijection
$J_n:X^n\to X^n$ such that $J_nS^{ii+1}_nJ_n^{-1}=(S')^{ii+1}_n$, where
$S^{ii+1}_n$ is the same as in Theorem \ref{twisted}. In this way, twisted
actions of $B_n$ given by $S$ and $S'$ are conjugate.
\end{theorem}
We remark that the statement (iv) of Theorem \ref{derived} was proved in
\cite{ESS} (cf. Prop. 1.7) for unitary solutions and in \cite{LYZ}
(cf. Th. 3) for injective solutions that are defined in Section 2.3.

{\bf \noindent Proof:}

It is easy to see that (iv) implies (ii). Statement (iii) follows from
definitions of structure group, derived structure group and $\phi$. Let us
show that (i) implies (iv).  Define $J_n$ inductively as $J_1=id_X$,
$J_{n}=Q_n(J_{n-1}\times id_X)$, where $Q_n:X^n\to X^n$ is defined
as $Q_n(x_1,...,x_n)=(x_n^{-1}*x_1,...,x_n^{-1}*x_{n-1},x_n)$. We prove
formula $J_nS^{ii+1}_nJ_n^{-1}=(S')^{ii+1}_n$ by induction on $n$. For
$n=2$ (induction base) the relation follows directly from definition of
$\phi$. Suppose the relation holds for $n=k$, let us prove it holds for
$n=k+1$. Since  $J_{k+1}=Q_{k+1}(J_{k}\times id_X)$ and $Q_{k+1}$ commutes
with $(S')^{ii+1}_{k+1}$ for $i<k$ (by (i)) the relation is true for
$i<k$. So, it remains to prove it for $i=k$ when it becomes identical to
the relation of the induction base. 
 
Now we have to show that (i) is true. It is enough to check
that for every 
$t,y,z\in X$, $t^{-1}*\phi(y,z)=\phi(t^{-1}*y,t^{-1}*z)$ or, equivalently
by (\ref{phi}):
\begin{equation}
\label{invariance}
t^{-1}*(z^{-1}*((y*z) \circ
y))=(t^{-1}*z)^{-1}*(((t^{-1}*y)*(t^{-1}*z))\circ (t^{-1}*y)).
\end{equation}
The linking relation of Theorem \ref{linking} states that for
$x, y, t\in X$
\begin{equation}
\label{rhs}
((y^{-1}*x)\circ t)^{-1}*(x\circ y)=((y\circ
t)^{-1}*x)\circ(t^{-1}*y)
\end{equation}
holds. If we substitute $x=y*z$ in relation (\ref{rhs}) we can
rewrite its right hand side  as follows:
$$
((y\circ t)^{-1}*(y*z))\circ (t^{-1}*y)=((y\circ
t)^{-1}*(yt*(t^{-1}*z)))\circ (t^{-1}*y)
$$
where the product $yt$ is the product of elements $y, t$ in group $G_X$,
thus $yt=(y\circ t)(t^{-1}*y)$ and
$$
((y\circ t)^{-1}*(yt*(t^{-1}*z)))\circ
(t^{-1}*y)=((t^{-1}*y)*(t^{-1}*z))\circ (t^{-1}*y).
$$
In this way, relation (\ref{rhs}) is equivalent to relation
(\ref{lhs}) below:
\begin{equation}  
\label{lhs}
(z\circ t)^{-1}*((y*z)\circ y)=((t^{-1}*y)*(t^{-1}*z))\circ (t^{-1}*y).
\end{equation}
Substituting, (\ref{lhs}) into the right side of (\ref{invariance}) we get
$$
(t^{-1}*z)^{-1}*(((t^{-1}*y)*(t^{-1}*z))\circ
(t^{-1}*y))=((t^{-1}*z)^{-1}(z\circ
t)^{-1})*((y*z)\circ y).
$$
Using the relation $t^{-1}z^{-1}=(t^{-1}*z)^{-1}(z\circ t)^{-1}$ in group
$G_X$ we verify that invariance relation (\ref{invariance}) is true.
The theorem is proved. $\blacksquare$

Note that by construction we have two natural maps
$\psi_G:X\to G_X$ and $\psi_A:X\to A_X$ that are not necessarily
embeddings (see Examples below). Theorem \ref{action} shows that
the latter map is $G_X$-invariant with respect to a suitable $G_X$-module
structure on $A_X$.  
Define the group $Permut(X)$ as the group of all
permutations of the set $X$. Then both $*$ and $\circ$ define
homomorphisms from $G_X$ to $Permut(X)$. Denote by $Aut_X(A_X)$ the group
of automorphisms of $A_X$ that map the generating set $\psi_A(X)$ onto
itself. 
\begin{theorem}
\label{action}
For a nondegenerate braided set $(X,S)$ the group homomorphism
$*:G_X\to Permut(X)$
can be uniquely lifted to the
homomorphism $\hat{*}:G_X\to Aut_X(A_X)$ such that for any $g\in G_X,\
x\in X$ $\psi_A(g*x)=g\hat{*}\psi_A(x)$. 
\end{theorem}
{\bf Proof:}

It is clear that the map $*:G_X\to Permut(X)$ can be uniquely extended
to the homomorphism from the group $G_X$ to the group of automorphisms of 
the free group generated by $X$. So, it
is enough to show that such an extension respects the generating relations
in group $A_X$. This immediately follows from statements (i), (iii) of
Theorem \ref{derived}. $\blacksquare$

{\bf Example}\cite{Dr}
Let $X$ be a union of conjugacy classes in a  group $G$,
i.e. $gXg^{-1}=X$ for any $g\in G$. Define $S:X\times X\to X\times X$ as
$S(x,y)=(xyx^{-1},x)$ for $x,y\in X$. It is easy to see that
$(X,S)$ is a braided nondegenerate
set. It was called the conjugate solution in \cite{LYZ}.
Clearly, this solution coincides with its own derived solution -
therefore $G_X=A_X$. The action of $G_X$ on $A_X$ is easily seen to be
trivial, i.e. $g\hat{*}h=h$ for any $g\in G_X,\ h\in A_X$. The embedding
$i:X\to G$ can be extended to a homomorphism $I:G_X\to G$ that maps
elements of $\psi_G(X)$ onto $X$, thus 
the map $\psi_G:X\to G_X$ is injective. So $(X, S)$ is an injective
solution (cf. Definition \ref{def6}).
 
\subsection{Bijective 1-cocycle}
Recall that we
started with a nondegenerate braided set $(X, S)$ and constructed two
groups $G_X$, $A_X$ and the action $\hat{*}$ of $G_X$ on $A_X$. In the
future we drop $\hat{}$ in $\hat{*}$ and denote the action of
$g\in G_X$ on $h\in A_X$ from Theorem \ref{action} simply by
$g*h$. Theorem \ref{cocycle} is the main step
towards classification of nondegenerate braided sets in group theoretical
terms.
\begin{definition}
Suppose a group $G$ acts on a group $A$ by automorphisms meaning that
there
is a homomorphism from $G$ to $Aut(A)$. Denote the product of elements
$g_1,\ g_2$ in $G$ by $g_1g_2$ and the product of elements $a_1,\ a_2$ in
$A$ by $a_1\bullet a_2$. We call a map $\pi:G\to A$ a
1-cocycle if for any $g_1, g_2\in G$ 
\begin{equation}
\pi(g_1g_2)=(g_2^{-1}*\pi(g_1))\bullet \pi(g_2) 
\end{equation}
where $*$ stands for the action of $G$ on $A$.
\end{definition} 

\begin{theorem}(cf. \cite{LYZ}, Th. 2 and \cite{ESS}, Prop. 2.5) 
\label{maintheorem}
For a nondegenerate braided set $(X, S)$ there exists a unique
bijective 1-cocycle $\pi:G_X\to A_X$ such that $\pi\psi_G=\psi_A$ on $X$.
\label{cocycle}
\end{theorem} 
\underline{Remark.} In \cite{LYZ} the authors introduced on $G_X$ another
group structure. It follows from our and their results that $G_X$ with the
new group structure is isomorphic to $A_X$ via $\pi:G_X\to A_X$.

{\bf \noindent Proof:}

{\bf Construction of the 1-cocycle.}

Let us construct the 1-cocycle $\pi:G_X\to A_X$. Consider the semidirect
product $G_X\ltimes A_X$. The group $G_X\ltimes A_X$ consists of pairs  
$(g, h)$, $g\in G_X,\ h\in A_X$ with the group operation given by the
formula
$$
(g, h)(g',h')=(gg',((g')^{-1}*h)\bullet h').
$$ 
Define the map $s:X\to G_X\ltimes A_X$ via the formula 
$s(x)=(\psi_G(x),\psi_A(x))$ for any $x\in X$. We claim that there exists
a group homomorphism $\bar{s}:G_X\to G_X\ltimes A_X$ such that
$\bar{s}\psi_G=s$. Indeed, since $\psi_G(X)$ generates the whole $G_X$ to
show that $\bar{s}$ exists it is necessary
and sufficient to check that $s$ respects the relations in $G_X$,
i.e. that for $x, y \in X$ $s(x)s(y)$ is equal to $s(x\circ
y)s(y^{-1}*x)$ in group $G_X\ltimes A_X$. When we
formally multiply terms out, the above condition transforms into the
relation
$$
(xy, (y^{-1}*x)\bullet y)=((x\circ y)(y^{-1}*x),((y^{-1}*x)^{-1}*(x\circ
y))\bullet (y^{-1}*x)),
$$
which coincides componentwise with the defining relations in groups $G_X$,
$A_X$. Let us define the projection $p:G_X\ltimes A_X\to A_X$ by the
formula $p(g,h)=h$. Introduce $\pi:G_X\to A_X$ as $\pi=p\bar{s}$. Clearly,
$\pi$ is a 1-cocycle. The tricky part is to show that $\pi$ is
bijective. 

{\bf Proof that $\pi$ is bijective.}
\begin{lemma} 
\label{t1t2}
(i) For $x_1,\ x_2\in X$, $g\in G_X$ if 
$\psi_G(x_1)=\psi_G(x_2)$ then $\psi_G(g*x_1)=\psi_G(g*x_2)$ and
$\psi_G(g\circ x_1)=\psi_G(g\circ x_2)$.

(ii)
For $x\in X$ $\psi_G((x*x)\circ x)=\psi_G(x*x)$ and 
$\psi_G((x^{-1}\circ x)^{-1}*x)=\psi_G(x^{-1}\circ x)$.

(iii)
For $x\in X$ $\psi_G(x)=\psi_G((x*x)^{-1}\circ (x*x))=\psi_G((x^{-1}\circ
x)*(x^{-1}\circ x))$.
\end{lemma}
To prove the lemma we notice that $S(x*x,x)=((x*x)\circ x,x)$ therefore 
$\psi_G((x*x)\circ x)\psi_G(x)=\psi_G(x*x)\psi_G(x)$ and
$\psi_G((x*x)\circ x)=\psi_G(x*x)$. Similarly 
$S(x,x^{-1}\circ x)=(x,(x^{-1}\circ x)^{-1}*x)$ implies that
$\psi_G((x^{-1}\circ x)^{-1}*x)=\psi_G(x^{-1}\circ x)$. So statement
(ii) of the lemma is proved. It is clear that statement (iii) follows from
(i) and (ii) thus it remains to prove (i). 

Let us show that
$\psi_G(x_1)=\psi_G(x_2)$ implies that $\psi_G(g*x_1)=\psi_G(g*x_2)$. For
a fixed $z\in X$ it
is enough to reason that $\psi_G(x_1)=\psi_G(x_2)$ if and only if
$\psi_G(z^{-1}*x_1)=\psi_G(z^{-1}*x_2)$. 
Since $S(x_1,z)=(x_1\circ z,z^{-1}*x_1)$ and 
$S(x_2,z)=(x_2\circ z,z^{-1}*x_2)$  we see that 

\begin{equation}
\label{1psi}
\psi_G(x_1)\psi_G(z)=\psi_G(x_1\circ z)\psi_G(z^{-1}*x_1)
\end{equation}
and
\begin{equation}
\label{2psi}
\psi_G(x_2)\psi_G(z)=\psi_G(x_2\circ z)\psi_G(z^{-1}*x_2).
\end{equation}
 Suppose
$\psi_G(x_1)=\psi_G(x_2)$  then since $*$ is an action of $G_X$ on $X$
we immediately see that $\psi_G(x_1\circ z)=\psi_G(x_2\circ z)$ therefore
equations (\ref{1psi})-(\ref{2psi}) imply that
$\psi_G(z^{-1}*x_1)=\psi_G(z^{-1}*x_2)$. Conversely, assume that
$\psi_G(z^{-1}*x_1)=\psi_G(z^{-1}*x_2)$. Our task is to prove that
$\psi_G(x_1)=\psi_G(x_2)$. We use the statement (iv) of Theorem
\ref{derived} in the following form: $S=J_2^{-1}S'J_2$, where 
$J_2(x,y)=(y^{-1}*x,y)$. We notice that 
$S'J_2(x_1,z)=(\phi(z,z^{-1}*x_1),z^{-1}*x_1)$ and
$S'J_2(x_2,z)=(\phi(z,z^{-1}*x_2),z^{-1}*x_2)$. Since $\phi(z,z^{-1}*x_1)$
is the action of $\psi_A(z^{-1}*x_1)\in A_X$ on $z\in X$ and
$\psi_A=\pi\psi_G$ we see that $\phi(z,z^{-1}*x_1)=\phi(z,z^{-1}*x_2)$. In
this way, first components of $S(x_1,z)=J_2^{-1}S'J_2(x_1,z)$ and
$S(x_2,z)=J_2^{-1}S'J_2(x_2,z)$ coincide, i.e. $x_1\circ z=x_2\circ z$. So
the equations (\ref{1psi})-(\ref{2psi}) imply that
$\psi_G(x_1)=\psi_G(x_2)$. In a similar fashion one can show that
$\psi_G(x_1)=\psi_G(x_2)$ if and only if $\psi_G(z\circ x_1)=\psi_G(z\circ
x_2)$. Lemma is proved. $\blacksquare$.

We aim to construct the map $h:A_X\to G_X$ inverse to $\pi$. At first we
define $h$ on $F_X$ - the free group generated by $X$ and then show that
it descends to $A_X$. We note that $G_X$ acts on $F_X$ via $f\to g*f$ for 
$f\in F_X,\ g\in G_X$ - the action induced from *-action of $G_X$ on $X$.
For $x\in X\subset F_X$ define
$$
h(x)=\psi_G(x),\ h(x^{-1})=(\psi_G(x^{-1}\circ x))^{-1}
$$
Notice that $x^{-1}\circ x$ is the same as $\psi_G(x)^{-1}\circ x$.
For an element $f=x_1\bullet ...\bullet x_k\in F_X$ of length $k$ define
inductively
$h(x_1\bullet ...\bullet x_k)=h(h(x_2\bullet
...\bullet x_k)*x_1)h(x_2\bullet ...\bullet x_k)$, where for each $i$, 
$x_i\in
X\subset F_X$ or $x_i^{-1}\in X\subset F_X$. In the above
$f=x_1\bullet ...\bullet x_k$
was the minimum decomposition of $f\in F_X$ and $k$, the length of f, is
the number of elements in such a decomposition. 
The only element of length 0
in $F_X$ is identity $e$ and we put $h(e)=1\in G_X$. We claim that for
$a,\ b\in F_X$
\begin{equation}
\label{h}
h(a\bullet b)=h(h(b)*a)h(b).
\end{equation}
Indeed, we can verify (\ref{h}) by induction on the length of
$a$. Suppose we
know that (\ref{h}) holds for for all elements $a$ of length $k$ and
want to check that $h(a\bullet y\bullet b)=h(h(b)*(a\bullet y))h(b)$ for
$y\in X\subset F_X$ or
$y^{-1}\in X\subset F_X$. We simplify $h(h(b)*(a\bullet y))h(b)$
as follows:
\begin{eqnarray*}
h(h(b)*(a\bullet y))h(b)=h((h(b)*a)\bullet (h(b)*y))h(b)\\
=h(h(h(b)*y)*(h(b)*a))h(h(b)*y)h(b). 
\end{eqnarray*}
On the other hand, we know that 
\begin{eqnarray*}
h(a\bullet y\bullet
b)=h(h(y\bullet b)*a)h(y\bullet b)=h(h(h(b)*y)h(b)*a)h(h(b)*y)h(b)\\
=h(h(h(b)*y)*(h(b)*a))h(h(b)*y)h(b).
\end{eqnarray*}
We see that $h(a\bullet y\bullet b)=h(h(b)*(a\bullet y))h(b)$. In this
way,
we just need to check
the induction base, namely that $h(y\bullet b)=h(h(b)*y)h(b)$. There are
two cases
to consider:
\begin{enumerate}
\item $Length(y\bullet b)=1+ Length(b)$ 
\item $Length(y\bullet b)=Length(b)-1$, i.e. $b=y^{-1}\bullet b'$ and
$Length(b)=Length(b')+1$
\end{enumerate}

In the first case the formula $h(y\bullet b)=h(h(b)*y)h(b)$ follows from 
definition of
$h$. In the second case without loss of generality assume that $y\in X$.
Then, $h(b)=h(y^{-1}\bullet b')=h((h(b')*y)^{-1})h(b')$ and letting
$h(b')*y=z$
we get
\begin{eqnarray*}
h(h(b)*y)h(b)=h(h((h(b')*y)^{-1})h(b')*y)h((h(b')*y)^{-1})h(b')\\
=h(h(z^{-1})*z))h(z^{-1})h(b')
=\psi_G(\psi_G(z^{-1}\circ z)^{-1}*z)\psi_G(z^{-1}\circ z)^{-1}h(b')\\
=(by\ Lemma)\ \psi_G((z^{-1}\circ z)\psi_G(z^{-1}\circ
z)^{-1}h(b')=h(b')=h(y\bullet b). 
\end{eqnarray*}

Now we are going to check that $h$ descends to $A_X$. Recall that $A_X$ is
given by generators - elements of $X$ and relations
$(z^{-1}*((y*z)\circ y))\bullet z=z\bullet y$ for each $y,\ z \in X$. 
We see that 
\begin{equation}
\label{pq}
h((z^{-1}*((y*z)\circ y))\bullet z)=\psi_G((y*z)\circ y)h(z)=h(z\bullet
y).
\end{equation}
If we define $p(y,z)=(z^{-1}*((y*z)\circ y))\bullet z\in F_X$ and
$q(y,z)=z\bullet y\in F_X$ then for any $g\in G_X$,
$h(g*(p(y,z)^{-1}\bullet q(y,z)))=h(g*(q(y,z)^{-1}\bullet p(y,z)))=1$.

Indeed, we
checked in the proof of Theorem \ref{action} that the defining relations
of group $A_X$ are $G_X$ -invariant with respect to * action, i.e. 
$g*p(y,z)=p(g*y,g*z)$ and $g*q(y,z)=q(g*y,g*z)$. Formula
(\ref{pq}) states that $h(p(y,z))=h(q(y,z))$, therefore 
\begin{eqnarray*}
h(p(y,z)^{-1}\bullet q(y,z))=h(h(q(y,z))*p(y,z)^{-1})h(q(y,z))\\
=h(h(p(y,z))*p(y,z)^{-1})h(p(y,z))=h(p(y,z)^{-1}\bullet p(y,z))=1.
\end{eqnarray*}
Thus
$h(g*(p(y,z)^{-1}\bullet
q(y,z)))=h(p(g*y,g*z)^{-1}\bullet q(g*y,g*z))=1$. Similarly
$h(g*(q(y,z)^{-1}\bullet p(y,z)))=1$.

In order to check that
$h:F_X\to G_X$ descends to the map $h:A_X\to G_X$ it is enough to check
that for any $a,b\in F_X,\ y,\ z\in X$ 
$h(a\bullet p(y,z)^{-1}\bullet
q(y,z)\bullet b)=h(a\bullet q^{-1}(y,z)\bullet p(y,z)\bullet
b)=h(a\bullet b)$ holds.  By formula
(\ref{h}) 
\begin{eqnarray*}
h(a\bullet 
p^{-1}(y,z)\bullet 
q(y,z)\bullet b)=h(h(p(y,z)^{-1}\bullet
q(y,z)\bullet b)*a)\\
h(p(y,z)^{-1}\bullet q(y,z)\bullet b), 
\end{eqnarray*}
so it is enough to show that $h(p^{-1}(y,z)\bullet
q(y,z)\bullet b)=h(b)$. But
$h(p^{-1}(y,z)\bullet q(y,z)\bullet
b)=h(h(b)*(p^{-1}(y,z)\bullet q(y,z)))h(b)=h(b)$ holds. 
Similar reasoning allows to check that
$h(a\bullet q^{-1}(y,z)\bullet p(y,z)\bullet b)=h(a\bullet b)$.
In this
way we constructed the map $h:A_X\to G_X$ such that
\begin{enumerate}
\item $h(a\bullet b)=h(h(b)*a)h(b)$
\item $h\psi_A=\psi_G$
\end{enumerate}
Condition 1 above implies that both $\pi h:A_X\to A_X$ and $h\pi:G_X\to
G_X$ are group homomorphisms while Condition 2 implies that these
homomorphisms are identities if restricted to corresponding generating
sets $\psi_A(X),\ \psi_G(X)$, thus $\pi h=id_{A_X},\ h \pi=id_{G_X}$.
Theorem \ref{maintheorem} is proved. $\blacksquare$

It turns out then that the groups $G_X$ and
$A_X$ both have abelian subgroups of finite index when $X$ is
finite. Namely, define 
$\Gamma=\{g\in G_X|g*x=x,\ g\circ x=x\ for\ all\ x\in X\}$. In other
words, $\Gamma$ is the intersection of the kernels of left and right
actions from Theorem \ref{linking}.
\begin{theorem} (cf. Section 2.5 in \cite{ESS} and Prop. 6 in \cite{LYZ})
\label{gamma}

(i) $\pi(\Gamma)$ is a normal $G_X$-invariant (w.r.t. to
$*$-action) subgroup lying in the center of $A_X$. $\Gamma$ is a normal
abelian subgroup in $G_X$, and $\pi:\Gamma\to \pi(\Gamma)$ is an
isomorphism.

(ii) The 1-cocycle $\pi:G_X\to A_X$ can be factored out through $\Gamma$
giving rise to the bijective 1-cocycle $\bar{\pi}:G_X/\Gamma \to
A_X/\pi(\Gamma)$.

(iii) If X is finite then both $G_X/\Gamma$ and $A_X/\pi(\Gamma)$ are
finite groups.
\end{theorem}
{\bf Proof:}

(i) Let us show that $\pi(\Gamma)$ is $G_X$-invariant and central
subgroup in $A_X$. From the defining relations in $G_X$ we see that
$\psi_G(x\circ y)=\psi_G(x)\psi_G(y)\psi_G(y^{-1}*x)^{-1}$ for all $x,
y\in X$. In fact, one can easily
check
that the above formula can be generalized for $x\in G_X$, $y\in X$.

\begin{lemma}
\label{gy}
For $g\in G_X, y\in X$
$$
\psi_G(g\circ y)=g\psi_G(y)\pi^{-1}(\psi_G(y)^{-1}*\pi(g))^{-1}.
$$
\end{lemma}
{\bf Proof of Lemma \ref{gy}.}
Let us show that if the statement of the Lemma holds for some $g\in G_X$
then it holds for $g^{-1}$. We need to check that 
\begin{equation}
\label{qqq}
g^{-1}\psi_G(g\circ y)\pi^{-1}(\psi_G(g\circ y)^{-1}*\pi(g^{-1}))^{-1}=y.
\end{equation}
Notice that $\pi(g^{-1})=(g* \pi(g))^{-1}$ therefore 
$$\psi_G(g\circ y)^{-1}*\pi(g^{-1})=((\psi_G(g\circ y)^{-1}g)*
\pi(g))^{-1}.$$ Since the statement of the Lemma holds for $g$, we have 
\begin{eqnarray*}
\pi^{-1}((\psi_G(g\circ y)^{-1}g)*
\pi(g))^{-1}=
\pi^{-1}((\pi^{-1}(\psi_G(y)^{-1}*\pi(g))\psi_G(y)^{-1})*\pi(g))^{-1}\\
=\pi^{-1}(\pi^{-1}(\psi_G(y)^{-1}*\pi(g))*(\psi_G(y)^{-1}*\pi(g)))^{-1}\\
=\pi^{-1}(\psi_G(y)^{-1}*\pi(g)).
\end{eqnarray*}
This implies the validity of equation (\ref{qqq}). Similar, a simple
computation shows that if the statement of the lemma is true for $g=g_1$
and $g=g_2$ then it is true for $g=g_1g_2$. In this way since by
definition of $G_X$ Lemma is
true for $g\in \psi_G(X)$ we prove the Lemma. $\blacksquare$  

In particular, for $g\in \Gamma,\ y\in \psi_G(X)$ one has 
$y=gy(\pi^{-1}(y^{-1}*\pi(g)))^{-1}$. Thus,
\begin{equation}
\label{center}
y^{-1}gy=\pi^{-1}(y^{-1}*\pi(g)).
\end{equation}
The condition ($\ref{center}$) implies that $\pi{\Gamma}$ belongs to the
center of $A_X$. Indeed, applying $\pi$ to relation
$y^{-1}g=\pi^{-1}(y^{-1}*\pi(g))y^{-1}$ we get that $\pi(y^{-1})$ commutes
with $\pi(g)$. Moreover, for $g\in \Gamma, h\in G_X$ 1-cocycle condition
implies that
\begin{equation}
\label{extension}
\pi^{-1}(\pi(h)\bullet \pi(g))=hg.
\end{equation}
So, in particular the product of elements in $\pi(\Gamma)$ is in
$\pi(\Gamma)$ and $\pi$ restricted to $\Gamma$ becomes an isomorphism
between
$\Gamma$ and $\pi(\Gamma)$. 

(ii) The relation (\ref{extension}) shows that $\pi$ can be lifted to
$\bar{\pi}:G_X/\Gamma\to A_X/\pi(\Gamma)$.

(iii) The kernels of each of the actions $*$, $\circ$ are of finite
indexes since the corresponding quotients are isomorphic to subgroups in
Permut(X). So, the intersection of kernels, i.e. subgroup $\Gamma$ has
finite index as well. $\blacksquare$ 

\begin{definition}
We call a 7-tuple 
$$(G, A, X, \rho_{GA}, \rho_{GAX},
\bar{\pi}, \overline{\psi_A})$$ a 
bijective cocycle 7-tuple if $G,\ A$ are groups, $\rho_{GA}:G\to Aut(A)$
is
an action of $G$ on $A$, $\rho_{GAX}:G\ltimes A\to Permut(X)$ is an action
of $G\ltimes A$ on $X$,
$\bar{\pi}:G\to A$ is a bijective 1-cocycle, $\overline{\psi_A}:X\to A$ is
$G\ltimes
A$-equivariant, where $G\ltimes A$ acts on $A$ by conjugation.
\end{definition}
Note that the
action of $G\ltimes A$ gives rise to two actions $\rho_{GX}:G\to
Permut(X)$
and $\rho_{AX}:A\to Permut(X)$.
Starting with a bijective cocycle 7-tuple, let $*:X\times X\to X$ and
$\phi:X\times
X\to X$ be defined as
$y*x=\rho_{GX}(\overline{\psi_G}(y))(x)$ and 
$\phi(x,y)=\rho_{AX}(\overline{\psi_A}(y))(x)$, where 
$\overline{\psi_G}=\bar{\pi}^{-1}\overline{\psi_A}$. 
Define 
$S:X\times X\to X\times X$ and $S':X\times X\to X\times X$
via $S(x,y)=(x\circ y, y^{-1}*x)$ and 
$S'(x,y)=(\phi(y,x), x)$, where 
$x\circ y$ is defined in such a way that
$\phi(y,x)=x^{-1}*((y*x)\circ y)$ holds. 

\begin{lemma}
\label{faith}
(i) For an arbitrary bijective cocycle 7-tuple, $(X, S)$ constructed above 
is a braided nondegenerate set.

(ii) $(X, S')$ is a derived nondegenerate braided set corresponding to 
$(X,S)$.
\end{lemma}
{\bf\noindent Proof:}

It is obvious from the definition of a bijective cocycle 7-tuple that $(X, 
S')$ is a braided nondegenerate set. We introduced $x\circ y$ in such a
way that $(X, S')$ becomes a derived solution corresponding to $(X, S)$
as long as we are able to show that $(X, S)$ is braided nondegenerate
itself. Since $\phi$ is $G$-invariant the argument identical to the one
made in the proof of Theorem \ref{derived} mitigates that for each integer
$n\geq 2$ there exists a bijection $J_n:X^n\to X^n$ such that 
$J_nS_n^{ii+1}J_n^{-1}=(S')_n^{ii+1}$. This implies that $S$ is
bijective and satisfies the braid relation (\ref{braid}). 

It remains to prove that $(X, S)$ is nondegenerate, i.e. that $x\circ y$
depends bijectively on $y\in X$ for a fixed $x\in X$. For that purpose we
show that we can define $g\circ y\in X$ for $y\in X,\ g\in G$
such that for any $x,\ y \in X$, $g, h\in G$ 
$x\circ y=\overline{\psi_G}(x)\circ y$ and  $(gh)\circ y=g\circ(h\circ
y)$,
i.e. $\circ$ is an action of $G$ on $X$.

Let a group homorphism $P:G\to G\ltimes A$ be
defined by the formula $P(g)=(g,\bar{\pi}(g))$ for $g\in G$. Then the map
$\rho^*_{GX}=\rho_{GAX}P$ is an action of
$G$ on $X$. For $g\in G, y\in X$ define 

$$
g\circ
y=\rho^*_{GX}(\bar{\pi}^{-1}(\rho_{GA}(\overline{\psi_G}(y)^{-1})(\bar{\pi}(g))))(y).
$$ 
We want to check $(gh)\circ y=g\circ(h\circ y)$. For notational convience
let
$\rho_{GA}(\overline{\psi_G}(y)^{-1})(\bar{\pi}(g))=\overline{\psi_G}(y)^{-1}*\bar{\pi}(g)$,
i.e. by $*$ we will mean both actions $\rho_{GA}$ and $\rho_{GX}$. Then
due to the fact that $P$ is a homomorphism it
is enough to verify that 
$$
\bar{\pi}^{-1}(\overline{\psi_G}(h\circ y)^{-1}*\bar{\pi}(g))=
\bar{\pi}^{-1}(\overline{\psi_G}(y)^{-1}*\bar{\pi}(gh))(\bar{\pi}^{-1}
(\overline{\psi_G}(y)^{-1}*\bar{\pi}(h))^{-1}
.$$  
Since $\bar{\pi}(gh)=(h^{-1}*\bar{\pi}(g))\bar{\pi}(h)$,
\begin{eqnarray*}
\bar{\pi}^{-1}(\overline{\psi_G}(y)^{-1}*\bar{\pi}(gh))=
\bar{\pi}^{-1}((\overline{\psi_G}(y)^{-1}h^{-1}*\bar{\pi}(g))(\overline{\psi_G}(y)^{-1}*\bar{\pi}(h)))\\
=\bar{\pi}^{-1}(\bar{\pi}^{-1}(\overline{\psi_G}(y)^{-1}*\bar{\pi}(h))
\overline{\psi_G}(y)^{-1}h^{-1}*\bar{\pi}(g))\bar{\pi}^{-1}(\overline{\psi_G}(y)^{-1}*\bar{\pi}(h)).
\end{eqnarray*}
In this way, if we check that 
\begin{equation}
\label{xyz}
\bar{\pi}^{-1}(\overline{\psi_G}(y)^{-1}*\bar{\pi}(h))\overline{\psi_G}(y)^{-1}h^{-1}=
\overline{\psi_G}(h\circ y)^{-1}
\end{equation}

we conclude that $g\circ (h\circ y)=(gh)\circ y$. Let us rewrite
(\ref{xyz}) as (cf. Lemma 2) 
$\overline{\psi_G}(h\circ
y)\bar{\pi}^{-1}(\overline{\psi_G}(y)^{-1}*\bar{\pi}(h))=h\overline{\psi_G}(y)$
and apply
$\bar{\pi}$ to it. We get that (\ref{xyz}) is equivalent to
$$
\overline{\psi_A}(\bar{\pi}^{-1}(\overline{\psi_G}(y)^{-1}*\bar{\pi}(h))^{-1}*(h\circ
y))(\overline{\psi_G}(y)^{-1}*\bar{\pi}(h))=(\overline{\psi_G}(y)^{-1}*\bar{\pi}(h))(\overline{\psi_A}(y)).
$$
The last equality follows from $A$-equivariance of $\overline{\psi_A}$
since
$\bar{\pi}^{-1}(\overline{\psi_G}(y)^{-1}*\bar{\pi}(h))^{-1}*(h\circ y)
= \rho_{AX}(\overline{\psi_G}(y)^{-1}*\bar{\pi}(h))(y)$ by definition of
$h\circ
y$. Lemma is proved. 

$\blacksquare$

One can combine two actions of $G$ on $X$ -  $\rho_{GX}$ and 
$\rho^*_{GX}$ 
into the action $\rho=\rho_{GX}\times \rho^*_{GX}:G\to Permut(X)\times
Permut(X)$ of $G$ on $X^2$.  
\begin{definition}
We call a bijective cocycle 7-tuple
$(G, A, X, \rho_{GA}, \rho_{GAX}, \bar{\pi}, \overline{\psi_A})$ faithful
if $\overline{\psi_G}(X)=\bar{\pi}^{-1}\overline{\psi_A}(X)$
generates $G$ and
the action $\rho: G\to Permut(X)\times Permut(X)$ is faithful.
\end{definition}

The following theorem is a characterization of finite braided
nondegenerate sets in group theoretical terms. We notice that this result
is a generalization of Proposition 2.11 in \cite{ESS}.

\begin{theorem}
\label{1-1}
The construction of Lemma \ref{faith} establishes a 1-1
correspondence between 
nondegenerate braided sets $(X, S)$ and  
faithful bijective cocycle 7-tuples.
\end{theorem}
{\bf\noindent Proof:}

Having a nondegenerate braided set $(X, S)$ it is straightforward to
construct the faithful bijective cocycle 7-tuple. Indeed, following the
notations of
Theorem \ref{gamma} we put $G=G_X/\Gamma$, $A=A_X/\pi(\Gamma)$ then 
$G$ acts faithfully on $X^2$ and we have a faithful cocycle
7-tuple. Conversely, if we have a faithful 7-tuple we can
construct a nondegenerate braided set by Lemma \ref{faith}. Let $G_X$ be
the structure group of so constructed nondegenerate braided set. It
follows from (\ref{xyz}) that there is a group homomorphism
$Q:G_X\to G$ such that $Q\psi_G=\overline{\psi_G}$. Since
$\overline{\psi_G}(X)$
generates
$G$, $Q$ is surjective, and since $G$ acts faithfully on $X^2$,
$Ker (Q)=\Gamma$, i.e. $G_X/\Gamma$ is isomorphic to $G$. $\blacksquare$ 
	
\subsection{Injective solutions}
In this section we talk about most tractable braided nondegenerate 
sets - injective solutions.

\begin{definition}
\label{def6}
We call a braided nondegenerate set  $(X,\ S)$ an injective solution if
the map $\psi_G:X\to G_X$ is injective.
\end{definition}

In \cite{LYZ} the authors noticed that in the absence of involutivity
the natural map $\psi_G:X\to G_X$ is not obviously injective which
creates difficulties in characterization of solutions. It turned
out that injectivity may indeed fail (see examples below). This motivates
Definition \ref{def6}.
\begin{lemma}
\label{inj}
A nondegenerate braided set $(X, S)$ is injective if and only if
its derived solution $(X,\ S')$ is injective.
\end{lemma}
{\bf Proof:}

The statement of Theorem \ref{maintheorem} implies that
injectivity of the map $\psi_A:X\to A_X$ is equivalent to injectivity
of $\psi_G:X\to G_X$. Since $A_X$ is the structure group of the derived
solution Lemma \ref{inj} is proved. $\blacksquare$

The importance of injective solutions is in the fact that their properties
and group-theoretical characterization are very similar to that of
involutive solutions \cite{ESS}. 
\begin{theorem}
\label{reverse}
(i) Let a group $G$ act on a group $A$ by $\rho_{GA}:G\to Aut(A)$
such that the bijective map 
$\pi:G\to A$ is a 1-cocycle. Then any $G\ltimes A$-invariant 
subset
$X \subset A$ has a
natural structure of a nondegenerate braided injective set given by
\begin{equation}
\label{sss}
S(x,y)=(\pi(\pi^{-1}(x)\pi^{-1}(y)(\pi^{-1}(\rho(\pi^{-1}(y)^{-1})(x)))^{-1}) ,\rho(\pi^{-1}(y)^{-1})(x)), 
\end{equation}
for $x,\ y \in X$. 

(ii)
Any nondegenerate braided injective set can be obtained by the
method
just described.
\end{theorem}
{\bf Proof:}

To prove (i) we can use Lemma \ref{faith}. Indeed,  we have the following
bijective cocycle 7-tuple: $G$, $A$, $X$, $\bar{\pi}=\pi$, $\rho_{GA}$ as
given; $\rho_{GAX}$ is induced from the adjoint action of $G\ltimes A$ on
its subgroup $A$, $\overline{\psi_A}=id_X$. It is straightforward to check
(cf. (\ref{xyz})) that the map $S$ constructed in Lemma \ref{faith}
coincides with the map $S$ given by formula (\ref{sss}). So it remains to
prove that the set $(X, S)$ is injective. Let $G_X$ be its structure
group. Arguing as in the proof of Theorem \ref{1-1} we see that there is a
group homomorphism $Q:G_X\to G$ such that $\pi Q\psi_G=id_X$. This implies
that $\psi_G$ is injective.

Conversely, if $(X,\ S)$ is  an injective  nondegenerate braided set
then $X=\psi_A(X)$ is a $G_X$-invariant (w.r.t. $*$ -
action) subset in $A_X$. Note that this subset is automatically
$G_X\ltimes A_X$-invariant. Recall that $S$ is given by the formula
$S(x,y)=(x\circ
y,y^{-1}*x)$ for $x,y\in X$ according to notations after Theorem
\ref{linking}. Let us make sure that the construction of Theorem
\ref{reverse}
yields the same map $S$ we already have. Indeed, for $x,y\in \psi_A(X)=X$, 
$\pi^{-1}(y)^{-1}*x=\psi_G(y)^{-1}*x=y^{-1}*x$ and
$\pi(\pi^{-1}(x)\pi^{-1}(y)(\pi^{-1}(\pi^{-1}(y)^{-1}*x))^{-1})=
\pi(\psi_G(x)\psi_G(y)\psi_G(y^{-1}*x)^{-1})=\pi(\psi_G(x\circ y))=
\psi_A(x\circ y)=x\circ y$. Theorem is proved. $\blacksquare$

Lemma \ref{inj} implies that injectivity of a given solution is
determined by the properties of the
function $\phi:X\times X\to  X$. In particular, for a symmetric set 
$(X, S)$, $\phi(y,x)=y$ and $A_X$ is the free abelian group generated by
$X$. Hence, symmetric sets are injective. We don't know any easy way to
check that a given function $\phi(y,x)$ corresponds to an injective
solution. While an injectivity criterion is provided by Theorem
\ref{criterion}, we give two simple necessary conditions 
below that are in many cases sufficient to check that a given solution is
not injective.
\begin{lemma}
\label{necessary}
If $(X, S)$ is an injective braided set then

(a) $\phi(x,x)=x$ for all $x\in X$,

(b) a pair $(y,x)\in X\times X$ satisfies $\phi(y,x)=y$ if and only if
$\phi(x,y)=x$.
\end{lemma}
Proof:

Suppose $(X,\ S)$ is injective. The group $A_X$ is generated by the
elements of $X$ subject to relations $\phi(y,x)\bullet y=y\bullet x$ for 
$x,\ y\in X$. Consequently 
$\psi_A(\phi(x,x))\bullet \psi_A(x)=\psi_A(x)\bullet \psi_A(x)$ and 
$\psi_A(\phi(x,x))=\psi_A(x)$ in
$A_X$. Since $\psi_A:X\to A_X$ is injective $\phi(x,x)=x$ on $X$. Now,
assume that $\phi(y,x)=y$. Then
$\psi_A(y)\bullet \psi_A(x)=\psi_A(x)\bullet \psi_A(y)$.
On the other hand, 
$\psi_A(\phi(x,y))\bullet \psi_A(y)=\psi_A(y)\bullet \psi_A(x)$, therefore
$\psi_A(\phi(x,y))=\psi_A(x)$ and $\phi(x,y)=x$. Lemma is
proved. $\blacksquare$   

{\bf Example.} Let $c,\ b\in Permut(X)$. Define $S:X\times X\to X\times X$
by the formula $S(x,y)=(by,cx)$. It is easy to see that $(X,\ S)$ is a
nondegenerate braided set if $bc=cb$. We claim that this solution is
injective if and only if $cb=id_X$. Indeed, suppose the solution is
injective, then by Lemma \ref{necessary} $\phi(x,x)=x$, i.e.
$x^{-1}*((x*x)\circ x)=cbx=x$. Therefore $cb=id_X$. Conversely, if
$cb=id_X$ then $(X,\ S)$ is symmetric and hence injective.

We remark that with each nondegenerate braided set we associated two
actions of the group $G_X$ ($*,\ \circ$) and an action of $A_X$ 
(via $\phi(x,y)$, which is a $\circ$ action of the derived solution) on
$X$. In particular, the latter action allows us to
construct a finite group $A_X^0\subset Permut(X)$, as the image of $A_X$
under that action, and a surjective homomorphism $p:A_X\to A_X^0$. Define
$M_X$ - a module over $A_X^0$ generated by $v_x$, $x\in X$ subject to
relations
$$
p(y)^{-1}v_x+v_y=p(x)^{-1}v_{\phi(y,x)}+v_x.
$$
By construction we have a natural map $\psi_M:X\to M_X$ given by $x\to
v_x$. It turns out that injectivity of $\psi_M$ is equivalent to
injectivity of $\psi_G$.
\begin{theorem}
\label{criterion}
(i) There exists a unique 1-cocycle $\theta:A_X\to M_X$ such that
$\theta\psi_A=\psi_M$, where $A_X$ acts on $M_X$ via $p$. 

(ii) $\theta$ is injective on $\psi_A(X)$.
\end{theorem}
Proof:

Statement (i) is clear from definitions of $A_X$ and $M_X$. Let us show
that (ii) holds. Let $Ker(p)=\Gamma_A$, i.e. $A_X^0=A_X/\Gamma_A$. 
Let $x_1, x_2\in \psi_A(X)$, $x_1\neq x_2$. We want to show that
$\theta(x_1)\neq \theta(x_2)$ in $M_X$. Fix a character 
$\xi:\Gamma_A\to \C^*$. Define a vector space 
$$V_\xi=\{f:A_X\to \C| f(a\bullet \gamma)=f(a)\xi(\gamma),\ \gamma\in
\Gamma_A\}.
$$
$V_\xi$ clearly has an $A_X$-module structure defined as
$(bf)(a)=f(b^{-1}\bullet a)$, where $f\in V_\xi,\ b,\ a\in A_X$. Choose a
lifting $g:A_X^0\to A_X$ (as a set only). Then $V_\xi$ is identified with
$Fun(A_X^0,\C)$, the space of functions on $A_X^0$  via 
$f\to f|_{g(A_X^0)}$. This space has a basis $\delta_a$, $a\in A_X^0$ such
that $\delta_a(b)=0$ if $a\neq b$ and $\delta_a(a)=1$. Let us define
$\varepsilon:\psi_A(X)\to Fun(A_X^0,\C)$ by the formula 
$x\delta_a=\delta_{p(x)\bullet a}\varepsilon(x)(a)$ for $x\in \psi_A(X)$,
$a\in A_X^0$.
We can always choose $\xi$ and $g$ such that $\varepsilon(x_1)\neq
\varepsilon(x_2)$. Indeed, if $p(x_1)=p(x_2)$ it
suffices to choose $\xi$ such that $\xi(x_1\bullet x_2^{-1})\neq 1$ for
any lifting $g$. If $p(x_1)\neq p(x_2)$ then for any character $\xi\neq 1$
there is a lifting $g$ such that $\varepsilon(x_1)\neq \varepsilon(x_2)$.
Moreover,
since $x\bullet y=\phi(y,x)\bullet x$ for $x,\ y\in \psi_A(X)$ and
$$
xy\delta_a=x\delta_{p(y)\bullet a}\varepsilon(y)(a)=
\delta_{p(x)\bullet p(y)\bullet a}\varepsilon(x)(p(y)\bullet
a)\varepsilon(y)(a),
$$
we get that 
$$
\varepsilon(x)(p(y)\bullet a)\varepsilon(y)(a)=
\varepsilon(\phi(y,x))(p(x)\bullet a)\varepsilon(x)(a).
$$ 
In this way, we have the equality 
$(p(y)^{-1}\varepsilon(x))\varepsilon(y)=
(p(x)^{-1}\varepsilon(\phi(y,x)))\varepsilon(x)$
in $Fun(A_X^0,\C)$. Hence we can construct an $A_X^0$-homomorphism 
from $M_X$ to $Fun(A_X^0,\C)$ given
by $v_z\to \varepsilon(\psi_A(z))$, $z\in X$. But 
$\varepsilon(x_1)\neq \varepsilon(x_2)$ thus $\theta(x_1)\neq
\theta(x_2)$.
\begin{corollary1}
The map $\psi_M:X\to M_X$ is injective if and only if $(X,\ S)$ is an
injective solution.
\end{corollary1}
\subsection{Rank of the structure group}

In this section we show how to compute the rank of the structure
group $G_X$ for a finite nondegenerate braided set $(X, S)$.

\begin{definition}

(i) The rank of a group G having an abelian subgroup of finite index
$\Gamma$ is defined as the rank of $\Gamma$.

(ii) We define the rank of a finite nondegenerate braided set $(X,\ S)$
to be the rank of its structure group $G_X$.   
\end{definition}
Clearly the above definition doesn't depend on the choice of $\Gamma$,
for any two abelian subgroups of finite index $\Gamma_1$, $\Gamma_2$
have the same rank that is equal to the rank of their intersection, which
has finite index in each of them. 
\begin{lemma}
The rank of a solution $(X,\ S)$ is equal to the rank of the derived
solution $(X,\ S')$.
\end{lemma}
Proof:

According to statement (i) of Theorem \ref{gamma} abelian subgroups 
$\Gamma$ and $\pi(\Gamma)$ of finite indexes in $G_X$ and $A_X$ are
isomorphic, thus $G_X$ and $A_X$ have the same rank. Lemma is proved.  

We aim to compute the rank of $A_X$.
Note that defining relations in group $A_X$ can be rewritten as
\begin{equation}
\label{xyx}
\phi(y,x)=x\bullet y \bullet x^{-1}.
\end{equation} 
Note that since $S(x,y)=(\phi(y,x),x)$ gives rise to a nondegenerate
braided set, $\phi(,x)$ can be extended to the action of $A_X$ on $X$. 
Introduce an equivalence
relation $\backsim$ on $X$ such that the orbits of the above action
become equivalence classes. Namely, consider the minimum equivalence
relation $\backsim$ such that $y\backsim \phi(y,x)$ for
$x,\ y \in X$. 

\begin{theorem}
\label{rank}
The rank of group $A_X$ is equal to the number of equivalence classes
with respect to equivalence relation $\backsim$ on $X$.
\end{theorem}

The proof relies on the following lemma. Let $H_1$, $H_2$, $H_3$ be three
groups, such that there is an exact sequence 
$$
1\to H_1\to H_2\to H_3\to 1
$$
and $H_1$ is cenral in $H_2$.
\begin{lemma}[\cite{CR}]
There exists an exact sequence (The  Hochschild - Serre
sequence) (\ref{serre}) for any
abelian
group B, where
$H^2(H_3,B)$ stands for the second cohomology group of $H_3$ with
coefficients in B. 
\end{lemma}
\begin{equation}
\label{serre}
1\to Hom(H_3,B)\to Hom(H_2,B)\to Hom(H_1,B)\to H^2(H_3,B)
\end{equation}
{\bf Proof of Theorem \ref{rank}.}

Let us use the sequence \ref{serre} in the following situation:
$H_2=A_X$,  $H_1\subset \pi(\Gamma)\subset A_X$ is the free abelian group
of rank $r$ of finite index in $\pi(\Gamma)$, $H_3=H_2/H_1$ - a finite
group, $B=\C^*$. Notice that $\pi(\Gamma)$ is central in $A_X$ therefore 
$H_1$ is central in $H_2$.
Then, since both $Hom(H_3,B)$ and $H^2(H_3,B)$ are finite
groups and $Hom(H_1,B)=(\C^*)^r$, the dimension of
$Hom(H_2,B)=Hom(A_X,\C^*)$ is equal
to $r$, the rank  of $A_X$. On the other hand, it is clear from
formula (\ref{xyx}) that  $Hom(A_X,\C^*)=(\C^*)^k$, where $k$ is the
number of equivalence classes in $X$. Theorem is proved. 

\begin{corollary1}
The rank of any solution $(X,\ S)$ is less or equal than $n$, the number
of elements in $X$, with the equality taking place if and only if $(X, S)$
is symmetric.
\end{corollary1}

{\noindent \bf Proof:}

It is clear that if rank of $A_X$ is equal to $n$ then equivalence
relation $\backsim$ on $X$ is trivial, $\phi(y,x)=y$ and the map
$\psi_A:X\to A_X$ is
injective. This immediately implies that the pair $(X, S)$ is symmetric.

{\bf Example.} Consider a permutation solution $S(x,y)=(bx,cy)$, 
$b,\ c\in Permut(X)$, $bc=cb$. It is easy to check that $\phi(y,x)=bcy$,
hence $y\backsim (bc)^my$ for any integer m. The rank of the permutation 
solution is equal to the number of equivalence
classes with respect to this relation, which, in turn,  is equal to the
number of independent cyclic permutations in canonical decomposition of
$bc$ (counting cyclic permutations of length 1).

\section{Linear And Affine Solutions}
\subsection{Linear Braided Sets}
In this section we
will look for nondegenerate braided sets of the following form:
$X$ is an abelian group, and $S$ is an
affine linear transformation of $X\times X$.
Such braided sets will be called {\it affine solutions}.
Considering affine solutions was motivated by the results in \cite{ESS}.

We will start with considering a special case, when
$S$ is an automorphism of $X\times X$. In this case, an affine solution
will be called {\it a linear solution}. For a linear solution,
$S$ has the form  
\begin{equation}
\label{linear}
S(x,y)=(ax+by,cx+dy),\ a,b,c,d\in\ End (X).
\end{equation}

It is easy to check that for $S$ of the form (\ref{linear})
the nondegeneracy is equivalent to invertibility of both $b$ and $c$ while 
braid relation is equivalent
to the equations \cite{Hi}

\begin{eqnarray}
\label{1st}
a(1-a)=bac,& d(1-d)=cdb,&\\
\label{2nd}
 ab=ba(1-d),&\ ca=(1-d)ac,&\ dc=cd(1-a),\\
\label{3rd}
bd=(1-a)db,&\  cb-bc=ada-dad.&
\end{eqnarray}

\begin{lemma}
\label{lemma}
Braided nondegenerate linear sets $(X,\ S)$
are in 1-1 correspondence with the quadruples
$(a,b,d,s)\in End(X)^4$ such that:

(i) $1-a,\ 1-d,\ b,\ 1+s$ are invertible,

(ii) s commutes with a,b,d and $sa=sd=0$,

(iii) $bdb^{-1}=(1-a)d,\ b^{-1}ab=a(1-d).$
The 1-1 correspondence is given via the formula
\end{lemma}
\begin{equation}
\label{s}
bc=(1-d+ad)(1-a)+s.
\end{equation}

{\bf \noindent Proof:}

Suppose $(a,b,c,d)$ solves (\ref{1st})-(\ref{3rd}). Note that the first of
equations (\ref{3rd}) implies that $bdb^{-1}=(1-a)d$, therefore
$b(1-d)b^{-1}=1-d+ad$. Moreover, if we
multiply the first equation of (\ref{2nd}) by $b^{-1}$ on the right we get
that $a=bab^{-1}b(1-d)b^{-1}$ thus relations (\ref{conj1}) hold. Similarly
from equations two and three of (\ref{2nd}) we obtain
relations (\ref{conj2}).
\begin{eqnarray}
\label{conj1}
bab^{-1}(1-d+ad)=a,\ bdb^{-1}=(1-a)d,\\
\label{conj2}
cac^{-1}=(1-d)a,\ cdc^{-1}(1-a+da)=d.
\end{eqnarray}

Above formulas show how to conjugate the elements of subalgebra generated
by $a$ and $d$ by elements $b$, $c$ and their products. 
We 
define $s$ from the relation (\ref{s}).
Notice that $sa=sd=0$. Indeed, multiplying (\ref{s}) by
$a$ on the right and using first of relations (\ref{1st}) and second of
relations (\ref{2nd}) we get that
$$sa=bca-(1-d+ad)(1-a)a=b(1-d)ac-(1-d+ad)bac.$$
Since according to (\ref{conj1}) $b(1-d)b^{-1}=1-d+ad$ we see that
$b(1-d)ac-(1-d+ad)bac=0$. Similarly $as=0$. The last of equations
(\ref{3rd}) imply that
\begin{equation}
\label{s2}
cb=(1-a+da)(1-d)+s.
\end{equation} 
Multiplying the relation (\ref{s2}) just obtained by $d$ we get that
$sd=ds=0$. Now we have everything to show that $1-a$, $1+s$ and $1-d$
are invertible. Indeed, $b,\ c$ are invertible because of nondegeneracy of
$(X,\ S)$. Since 
$$
bc=(1-d+s+ad)(1-a)=(1-a)(1-d+s+da)
,$$
$$
bc=(1-d+ad)(1-a+s)=(1-a+s)(1-d+ad)
,$$ 
and
$$
cb=(1-a+s+da)(1-d)=(1-d)(1-a+s+ad)
$$ 
we conclude that $1-a$ , $1-a+s$
and $1-d$
have right and left inverses and hence invertible. But 
$1-a+s=(1-a)(1+s)$, thus $1+s$ is invertible. 
Let us show that $s$ commutes with $c$. We
conjugate relation (\ref{s}) by $c$ and use relations (\ref{conj2}) to
conclude that 
\begin{equation}
\label{s3}
cb=(1-a+da)(1-d)+csc^{-1}.
\end{equation}
Comparing (\ref{s2}) to (\ref{s3}) we get that $s=csc^{-1}$, i.e. $s$
commutes with
$c$. Now, relation (\ref{s}) implies that 
$s$ commutes with $b$ as well. 
In this way,
starting from $(a,b,c,d)$, a solution to (\ref{1st})-(\ref{3rd}) we
constructed a quadruple $(a,b,d,s)$ satisfying conditions of the
lemma. Conversely, let us assume we start with $(a,b,d,s)$. Define
$c$ from the formula (\ref{s}). Since $1-d+ad=b(1-d)b^{-1}$ is invertible,
$c$
is
invertible as well.
It is straightforward to check
that so defined $(a,b,c,d)$ satisfy the braid relations
(\ref{1st})-(\ref{3rd}), and that $S:X\times X\to X\times X$ is bijective.
The lemma is proved.
\subsection{Injective Linear Solutions}
It turns out that injective linear solutions are easy to characterize.

\begin{theorem}
\label{injectivity}
A linear nondegenerate braided set of the form (\ref{linear}) is injective
iff $bc=(1-d+ad)(1-a)$ or, in the language of Lemma \ref{lemma}, $s=0$.
\end{theorem}
Proof:

Assume that $(X, S)$ is injective. Then, by Lemma \ref{necessary}
$\phi(x,x)=x$ on $X$. In the linear case it is easy to compute $\phi(y,z)$
explicitly:
$$\phi(y,z)=c((y*z)\circ y)+dz=c(a(y*z)+by)+dz=cac^{-1}(z-dy)+
cby+dz.$$
According to relation (\ref{s2}),
$cb=(1-a+da)(1-d)+s$. Plugging
it into formula for $\phi(y,z)$ and using (\ref{conj2}) we conclude that
$$\phi(y,z)=(1-(1-d)(1-a))z+((1-d)(1-a)+s)y.$$ Condition $\phi(x,x)=x$
immediately implies that $s=0$.

Conversely, assume that $s=bc-(1-d+ad)(1-a)=0$. Let us show that 
$\psi_A:X\to A_X$ is injective. The group $A_X$ is generated by elements
of $X$ subject to relations
\begin{equation}
\label{2.11}
\phi(y,z)\bullet z = z\bullet y,\ where\ y,z\in X\ and\ 
\phi(y,z)=z^{-1}*((y*z)\circ y).  
\end{equation} 
Denoting $K=(1-d)(1-a)$ we
see that the group $A_X$ is given by relations 
$((1-K)z+Ky)\bullet z=z\bullet y$. By Lemma \ref{lemma}, $K$ is
invertible. Let us now define the action of $A_X$
on $X$. We let the elements of generating set $\psi_A(X)$ act with $K$ on
$X$ and then extend this action to arbitrary elements of $A_X$. Consider a
semidirect product $A_X\ltimes X$ with respect to this action, and define
an embedding $J:X\to A_X\ltimes X$ given by the formula
$J(x)=(\psi_A(x),x)$. We notice that
$J((1-K)z+Ky)J(z)=J(z)J(y)$ in $A_X\ltimes X$. Indeed, 
\begin{eqnarray}\nonumber
J((1-K)z+Ky)J(z)=(\psi_A((1-K)z+Ky),(1-K)z+Ky)(\psi_A(z),z)\\
\nonumber =(\psi_A((1-K)z+Ky)\bullet \psi_A(z),K^{-1}((1-K)z+Ky)+z)\\
\nonumber =(\psi_A(z)\bullet \psi_A(y),K^{-1}z+y)=J(z)J(y).
\end{eqnarray}
So, $J$ can be extended to a homomorphism $\hat{J}:A_X\to A_X\ltimes
X$ such that $Pr\hat{J}\psi_A=id$, where $Pr:A_X\ltimes X\to
X$ stands for the projection to the second component. Therefore, the map
$\psi_A:X\to A_X$ is injective. Theorem is proved. 
\begin{corollary1}
(i)
Let $(X,\ S)$ be a nondegenerate braided linear set of the form
(\ref{linear}). Then the pair $(X,\ \hat{S})$ with 
$\hat{S}:X\times X\to X\times X$ given by
$\hat{S}(x,y)=(ax+by,cx+(d-s)y)$, where $s$ is defined in \ref{s},
is an injective solution.

(ii) Suppose $(X,\ S)$ is an injective linear solution and $s:X\to X$ 
satisfies $sa=as=0$, $sb=bs$, $sd=ds=-s^2$,
$sc=cs$. Then, $(X,\ \breve{S})$ with $\breve{S}:X\times X\to X\times X$
given by
$\breve{S}(x,y)=(ax+by,cx+(d+s)y)$ is a nondegenerate braided set that
corresponds under the correspondence of Lemma \ref{lemma} to the quadruple
$(a,b,d,s)$.
\end{corollary1}
{\bf \noindent Proof:}

(i) It is easy to check $(X,\ \hat{S})$ is a braided nondegenerate set by
directly checking relations
(\ref{1st}) - (\ref{3rd}). Since $s$ commutes with everything and
$sa=sd=0$ the above task is pretty simple. Also, it is obvious that
$bc=(1-(d-s)+a(d-s))(1-a)$, thus by Theorem \ref{injectivity} 
$(X,\ \hat{S})$ is injective. Proof of (ii) is similar to (i) and is left
to the reader.

{\bf \noindent Examples.} 

1. Consider the linear solution $S(x,y)=(cy,\ bx)$ with $c,\
b$ being linear automorphisms of $X$ subject to $cb=bc$. We see that 
s=bc-1, therefore $(X, \hat{S})$ with $\hat{S}(x,y)=(cy,bx+(1-bc)y)$ is an
injective solution. We obtain the same solution by a different method in
Example 2 at the end of this section. 

2. It was shown in \cite{ESS} that symmetric nondegenerate linear
solutions of the form \ref{linear} are given as the solutions to the
following equations:
$bab^{-1}=\frac{a}{a+1}$, $c=b^{-1}(1-a^2)$, $d=\frac{a}{a-1}$ with
$b,\ c$ being invertible. In
particular a large class of solutions of this kind considered corresponded
to
nilpotent $a$, i.e. there was $n$ such that $a^n=0$. Define
$s=ma^{n-1}$, $m$ being any integer. It is easy to see that  so defined
$s$ satisfies conditions of
Corollary 1 (ii) hence $S''(x,y)=(ax+by,cx+(d+ma^{n-1})y)$ is a
nondegenerate braided set, which is not symmetric and not injective unless
$ma^{n-1}=0$. 

\begin{theorem}
\label{affine}

Linear braided nondegenerate injective sets $(X, S)$ on an abelian group
$X$ 
are in 1-1 correspondence with triples
$(a,b,d)$ of endomorphisms of $X$ such that $b,\ 1-a,\ 1-d$ are invertible
and $bdb^{-1}=(1-a)d,\ b^{-1}ab=a(1-d)$.
\end{theorem}
Proof: Straightforward application of Lemma \ref{lemma} and Theorem
\ref{injectivity}.

\begin{corollary1}
Linear braided nondegenerate injective sets $(X, S)$ on an abelian group
$X$ are in 1-1 correspondence with triples
$(p,q,z)$ of automorphisms of $X$ such that $pq=qp$ and $z^2-z(p+q)+pq=0$.
\end{corollary1}

The statement of Corollary 1 follows from Theorem \ref{affine} via
change of variables $p=b^{-1},\ q=(1-a)(1-d)b^{-1},\ z=(1-a)b^{-1}$.

{\bf \noindent Examples.}

1. If 
$(X,S)$ is unitary then $(1-a)(1-d)=1$ and therefore $p=q$. In this way,
nondegenerate unitary linear braided sets are characterized as
representations of algebra
generated by invertible $p,\ z$ subject to $z^2-2zp+p^2=0$.

2. Put $z=p$, then $q$ can be anything as long as it is invertible and
$pq=qp$. Correspondingly, $a=0,\ b=p^{-1},\ 1-d=p^{-1}q,\ c=q$. The map
$S$ is defined in the following way: $S(x,y)=(p^{-1}y,qx+(1-p^{-1}q)y)$.
Similarly if we let $z=q$ we obtain the solution given by the formula:
$S(x,y)=(py+(1-q^{-1}p)x,q^{-1}x)$.

3. Let $\epsilon_1$ and $\epsilon_2$ be two nilpotent operators on $X$ of
order two, i.e. $\epsilon_1^2=\epsilon_2^2=0$. Then we let
$p=1+\epsilon_1$, $q=1-\epsilon_1$, $z=1+\epsilon_2$. It is easy to check
that so defined $p$ and $q$ commute and that $z^2-z(p+q)+pq=0$. We recover
$a$, $b$, $c$, $d$ from $\epsilon_1$ and $\epsilon_2$ via the formulas
$b=1-\epsilon_1$, $a=\epsilon_1-\epsilon_2+\epsilon_2\epsilon_1$,
$d=1-(1+\epsilon_1)(1-\epsilon_2)(1-2\epsilon_1)$,
$c=(1+\epsilon_1)(1-\epsilon_2)(1-\epsilon_1)(1+\epsilon_2)(1-\epsilon_1)$.
In particular, if $\epsilon_1\epsilon_2=\epsilon_2\epsilon_1$ the
corresponding linear solution has the form: 
$$
S(x,y)=((\epsilon_1-\epsilon_2+\epsilon_2\epsilon_1)x+(1-\epsilon_1)y,
(1-\epsilon_1)x+(\epsilon_1+\epsilon_2-\epsilon_2\epsilon_1)y)
$$
4.Suppose $X$ is a n-dimensional vector space and 
$p,\ q$ are invertible operators on it having $2n$ distinct
eigenvalues. Let $(v_1...v_n)$ be the basis of $X$ in which both $p$ and
$q$ are diagonalized, i.e. $pv_i=p_iv_i$, $qv_i=q_iv_i$. Then
$zv_i=z_iv_i$,
where for each $i$ either $z_i=p_i$ or $z_i=q_i$. Indeed, 
$0=(z^2-z(p+q)+pq)v_i=(z-p_i)(z-q_i)v_i$. The vector subspace of $X$
generated by applying $z$ to $v_i$ ($i$ is fixed) is annihilated by 
$(z-p_i)(z-q_i)$ therefore it has a basis of $p_i$ and $q_i$ eigenvectors
of $z$. Therefore, it has to be one dimensional, i.e. $zv_i=p_iv_i$ or
$zv_i=q_iv_i$. It is also easy to see that $z$ given by diagonal matrix
in $(v_1...v_n)$ basis with $p_i$ or $q_i$ entries on the $i$-th place
does satisfy $z^2-z(p+q)+pq=0$ on $X$.  
\subsection{Affine solutions}
In this section we talk about general affine solutions on an abelian group
$X$.

\begin{definition}
The solution (braided, nondegenerate) $(X,\ S)$ of the form
(\ref{genaffine}) is called affine.
\end{definition}
\begin{equation}
\label{genaffine}
S(x,y)=(ax+by+z,cx+dy+t),\ \ \ a,\ b,\ c,\ d\in End(X),\ t,\ z\in X.
\end{equation}
\begin{lemma}
\label{associated}
The pair $(X, S)$ of the form (\ref{genaffine}) is a solution if and only
if (\ref{1st})-(\ref{3rd}) and (\ref{4th}) hold. Therefore, any affine
solution gives rise to a linear solution $(X, S^*)$
with $S^*(x,y)=(ax+by,cx+dy)$.
\begin{equation}
\label{4th}
cdz+dt=0,\ az+bat=0,\ (c+d-ad-1)z+(da+1-a-b)t=0.
\end{equation}
\end{lemma}
Proof: Straightforward.
\begin{definition}
We call $(X, S^*)$ from Lemma \ref{associated} the linear part
of an affine solution $(X, S)$. Conversely, we call $(X, S)$ an affine
solution associated with $(X, S^*)$.
\end{definition}
\begin{theorem}
\label{aff}
(i) Let $(X,\ S^*)$ be a linear braided nondegenerate set,
$S^*(x,y)=(ax+by,cx+dy)$. Then,
$(X,\ S)$ given by (\ref{genaffine}) is an affine solution associated with
$(X,\ S^*)$ if and only if $t=-c(1-a)^{-1}z+k$ and
$ak=dk=0$, $(b-1)k = sz$, where s is defined in (\ref{s}).

(ii) An affine solution $(X,\ S)$ is injective if and only if its linear
part is injective and $k=0$ in the above characterization. In this way,
injective affine solutions associated with a given injective linear
solution $(X,\ S^*)$ are in 1-1 correspondence with elements $z\in X$, $t$
being given by  $t=-c(1-a)^{-1}z$.
\end{theorem}
Proof:

Let us proof part (i). Note that since
$cdz+dt=d(c(1-a)^{-1}z+t)$ and az+bat=$a(z+(1-a)c^{-1}t)$ by equations
(\ref{1st})-(\ref{2nd})
we can rewrite two of the relations (\ref{4th}) as 
$$
d(c(1-a)^{-1}z+t)=0,\ a(z+(1-a)c^{-1}t)=0.
$$
Therefore, if we define $k$ from the relation
\begin{equation}
\label{k}
t=-c(1-a)^{-1}z+k,
\end{equation}
we  see using (\ref{1st}) that $ak=dk=0$. Now, (\ref{conj2}) implies
$c(1-a)^{-1}c^{-1}=(1-a+da)^{-1}$ thus we can transform (\ref{k})
into $(1-a+da)t=-cz+(1-a+da)k=-cz+k$. In this way, we rewrite last of the
relations (\ref{4th}) as
\begin{equation}
\label{bt}
(d-ad-1)z-bt+k=0
\end{equation} 
If we substitute $t$ from (\ref{k}) into  (\ref{bt}) and use (\ref{s})
we get that $(b-1)k=sz$. Part (i) is proved. In order to prove part (ii)
we compute the function $\phi:X\times X\to X$. Define by $\phi^*$ the
corresponding function $\phi$ of the linear part $(X,\ S^*)$. Then, it is
easy to check that $\phi(y,x)=\phi^*(y,x)+k$. Now, assume that $(X, S)$ is
injective. Then by Lemma \ref{necessary} $\phi(x,x)=x$, hence
$\phi^*(x,x)=x$ and $k=0$. As we saw in the proof of Theorem
\ref{injectivity} $\phi^*(x,x)=x$ implies that $s=0$ and thus $(X, S^*)$
is
injective. Conversely, if $(X, S^*)$ is injective and $k=0$ then
$\phi(y,x)=\phi^*(y,x)$ and hence the derived solution of 
$(X, S)$ coincides with the derived solution of $(X, S^*)$ and hence is
injective. Theorem is proved.


\begin{thebibliography}{99}
\bibitem[Dr]{Dr} Drinfeld, V. Some unsolved problems in quantum group
theory, {\it Lecture Notes in Math, 1510}, p.1-8.
\bibitem[ESS]{ESS} Etingof, P.; Schedler, T.; Soloviev, A. Set-theoretical
solutions to the quantum Yang-Baxter equation. q-alg/9801047,
Duke Math. J. 1999.
\bibitem[CR]{CR} Methods of representation theory. Vol. I. With
applications to finite
groups and orders. Pure and Applied Mathematics. A Wiley-Interscience
Publication. John Wiley \& Sons, Inc., New York, 1981
\bibitem[Hi]{Hi} Hietarinta, J. Permutation-type solutions to the
Yang-Baxterand other simplex equations. Q-alg 9702006 1997.

\bibitem[LYZ]{LYZ} Lu, J-H.; Yan, M.; Zhu, Y-C. On set-theoretical
Yang-Baxter equation. To appear in Duke Math. J. 1999.
\end{thebibliography}
\end{document}